\documentclass[a4paper]{amsart}

\usepackage{times,amssymb,a4wide}

\swapnumbers

\newtheorem{thm}{Theorem}[section]
\newtheorem{lem}[thm]{Lemma}
\newtheorem{cor}[thm]{Corollary}
\theoremstyle{definition}

\numberwithin{equation}{section}

\newcommand{\eps}{\varepsilon}
\newcommand{\N}{\mathbb{N}}
\newcommand{\R}{\mathbb{R}}

\newcommand{\ip}[1]{\left<#1\right>}

\newcommand{\e}{\mathrm{e}}
\newcommand{\ds}{\,\mathrm{d}s}
\newcommand{\dt}{\,\mathrm{d}t}
\newcommand{\du}{\,\mathrm{d}u}

\begin{document}

\title[Iterates of Volterra Operators]%
{Asymptotic Behaviour of Iterates of Volterra Operators on $L^p(0,1)$}
\date{September, 2004}

\author{S.~P.~Eveson}
\address{Department of Mathematics\\University of York\\Heslington\\York
  YO10~5DD\\England}
\email{spe1@york.ac.uk}

\begin{abstract}
  Given $k\in L^1(0,1)$ satisfying certain smoothness and growth conditions at 
  $0$, we consider the Volterra convolution operator $V_k$ defined 
  on $L^p(0,1)$ by
  $$(V_ku)(t)=\int_0^t k(t-s)u(s)\ds,$$
  and its iterates $(V_k^n)_{n\in\N}$. We construct some much simpler 
  sequences which, as $n\to\infty$, are asymptotically equal in the operator 
  norm to $V_k^n$. This leads to a simple asymptotic formula for $\|V_k^n\|$ 
  and to a simple `asymptotically extremal sequence'; that is, a sequence 
  $(u_n)$ in $L^p(0,1)$ with $\|u_n\|_p=1$ and $\|V_k^nu_n\|\sim\|V_k^n\|$ as 
  $n\to\infty$. As an application, we derive a limit theorem for large
  deviations, which appears to be beyond the established theory.  
\end{abstract}

\maketitle

\section{Introduction}

A number of authors have recently published results on the asymptotic
behaviour of iterated Volterra operators on $L^2(0,1)$. Lao and
Whitley~\cite{LW97} established a number of estimates and provided
numerical  evidence for a conjecture about the operator norm of the
Riemann-Liouville  fractional integration operator which was
subsequently proved by  Kershaw~\cite{Kershaw99} and by Little and
Reade~\cite{LR98}. A somewhat stronger result was also independently
established by Thorpe~\cite{Thorpe98}. These results were generalised by
the author to other Volterra convolution operators on $L^2(0,1)$ and to
some extent to other Schatten-von Neumann norms in~\cite{Eveson03}.

We show here that analogues of most of these $L^2$ results also 
hold in $L^p$. The main result, Theorem~\ref{thm:main}, is that if
$k(t)=t^rf(t)$ where $r>-1$ and $f$ is differentiable at $0$, and we define
$$(V_ku)(t)=\int_0^t k(t-s)u(s)\ds,$$
then the asymptotic behaviour of $V_k^n$ is the same as that of $V_h^n$
where 
$$h(t)=f(0)t^r\e^{(k'(0)/k(0))t}.$$ 
The significance of this kernel is that there is a simple formula for its 
convolution powers, which leads to another asymptotically equivalent sequence
of operators of rank $1$ (Corollary~\ref{cor:power-exponential}), and an 
asymptoptic formula for the operator norm:
$$\|V_k^n\|\sim\frac{C_p(|f(0)|\Gamma(r+1))^n\e^{f'(0)/f(0)}}{\Gamma((r+1)n+1)}$$
as $n\to\infty$, where $C_p$ is a constant depending only on $p$, defined below.

As an application of these results, we derive a limit theorem for large
deviations (Section~\ref{sec:large}). The exact asymptotic formula for
$V_k^n$ may also have other applications: for example, it has recently
been shown~\cite{GMpre} that the Volterra operator $V$ with kernel $1$ is not 
supercyclic on any $L^p$ space; since the proof depends on direct calculations 
on the iterates $V^n$, the results and techniques established below
might lead to more general results on the same lines.

\section{Notation}
\label{sec:notation}

The term `sequence' will be applied equally to sequences indexed by
natural numbers or to generalised sequences indexed by positive real numbers.

Throughout, $p$ will denote a real number in the range $[1,\infty]$ and $q$
its H\"older conjugate, so $1/p+1/q=1$ if $1<p<\infty$ and $1$ is
conjugate to $\infty$. We use $\|\cdot\|_p$ to denote the norm on $L^p$
and the norm in the algebra of bounded operators acting on $L^p$. The duality 
pairing between $L^p(0,1)$ on $L^q(0,1)$ will be written using angle brackets: 
$\ip{f,g}=\int_0^1fg$. We denote by $C_p$ the constant
$$C_p=\begin{cases}
        \frac{1}{p^{1/p}q^{1/q}} & \text{ if } 1<p<\infty \\
        1                        & \text{ if } p=1 \text{ or } p=\infty.
      \end{cases}$$

Convolution of suitable functions on $[0,1]$ is defined by
$$(f*g)(t)=\int_0^tf(t-s)g(s)\ds$$
for $t\in[0,1]$ and the $n$-fold convolution power of $f$ 
is denoted by $f^{*n}$. For $k\in L^1(0,1)$, the Volterra convolution operator
$V_k$ associated with $k$ is defined on $L^p(0,1)$ by $V_kf=k*f$; it is well 
known that for any $p$, $V_k$ is a bounded operator on $L^p(0,1)$ with 
operator norm $\|V_k\|_p\leq\|k\|_1$.

If $(a_n)$ and $(b_n)$ are sequences of numbers, we shall say as usual
that $(a_n)$ and $(b_n)$ are asympotically equal, written $a_n\sim b_n$
as $n\to\infty$, if $a_n/b_n\to 1$ as $n\to\infty$.

Extending this idea to vectors, if $(u_n)$ and $(v_n)$ are sequences in
a normed linear space, we shall say that $u_n\sim v_n$ as $n\to\infty$ if
$$\frac{\|u_n-v_n\|}{\|u_n\|}\to 0.$$
It is easy to check that this is an equivalence relation. If $(T_n)$ is a 
sequence of bounded operators on a normed linear space and $(u_n)$ a sequence 
of non-zero vectors, we shall call $(u_n)$ \emph{asymptotically extremal} for
$(T_n)$ if $\|T_nu_n\|\sim\|T_n\|\|u_n\|$ as $n\to\infty$. We shall make 
frequent use of the following simple facts:
\begin{lem}\label{lem:parallel}
  \mbox{}
  \begin{enumerate}
  \item If $(u_n)$ and $(v_n)$ are sequences in a normed space $X$ and
    $u_n\sim v_n$ as $n\to\infty$, then $\|u_n\|\sim\|v_n\|$ as $n\to\infty$;
  \item if in addition $(S_n)$ and $(T_n)$ are sequences of bounded linear 
    operators on $X$ such that $S_n\sim T_n$ as $n\to\infty$ and $(u_n)$ is 
    asymptotically extremal for $(S_n)$ then $(v_n)$ is asymptotically 
    extremal for $(T_n)$.
  \end{enumerate}
\end{lem}
For sequences of positive real numbers we also use the notation
$a_n\lesssim b_n$ as $n\to\infty$ to mean that 
$\limsup_{n\to\infty}a_n/b_n\leq 1$, and $a_n\gtrsim b_n$ as $n\to\infty$
to mean that $\liminf_{n\to\infty}a_n/b_n\geq 1$.

\section{Kernels of the form $t^r\e^{\mu t}$}
\label{sec:power-exponential}

It is easy to check using the Laplace transform that if $k(t)=t^r\e^{\mu t}$
for some $r,\mu\in\R$ with $r>-1$, then the $n$-fold convolution power of
$k$ is given by
$$k^{*n}(t)=\frac{(\Gamma(r+1))^n}{\Gamma((r+1)n)}t^{(r+1)n-1}\e^{\mu t}$$
and we can choose to make this the definition of $k^{*n}$ for
non-integer $n>0$. For such kernels we can approximate $k^{*n}$ by
operators of rank $1$ and thus obtain asymptotic results. In fact, we need only
consider $k_0(t)=\e^{\mu t}$, because
$$k_0^{*n}(t)=\frac{1}{\Gamma(n)}t^{n-1}\e^{\mu t}$$
so $k^{*n}=(\Gamma(r+1))^nk_0^{*(r+1)n}$.

Throughout this section, $S_\lambda$ and $T_\lambda$ denote the operators on
$L^p(0,1)$ defined for any $\lambda\in\R$ by
\begin{align*}
  (S_\lambda f)(t)=\int_0^1\e^{\lambda(t-s)}f(s)\ds \\
  (T_\lambda f)(t)=\int_0^t\e^{\lambda(t-s)}f(s)\ds.
\end{align*}
We also write $e_\lambda$ for the function $t\mapsto\e^{\lambda t}$.
\begin{lem}\label{lem:rank-1}
  For any $p\in[1,\infty]$,
  $$\|S_\lambda\|_p\sim C_p\frac{\e^\lambda}{\lambda}$$
  as $\lambda\to\infty$ through $\R^+$. (The constant $C_p$ is defined in 
  Section~\ref{sec:notation}.) If $f_\lambda$ is defined by
  $$f_\lambda(t)=\begin{cases}
                  \e^{-g(\lambda)\lambda t} & \text{ if } p=1 \\
                  \e^{-\lambda t/(p-1)}     & \text{ if } 1<p<\infty \\
                  1                         & \text{ if } p=\infty
                 \end{cases}$$
  where $g$ is any function such that $g(\lambda)\to\infty$ as 
  $\lambda\to\infty$, then $(f_\lambda)$ is asymptotically extremal for
  $(S_\lambda)$.
\end{lem}
\begin{proof}
  We can write $S_\lambda$ in the form
  $$(S_\lambda f)(t)=\e^{\lambda t}\int_0^1\e^{-\lambda s}f(s)\ds
    =\ip{f,e_{-\lambda}}e_\lambda$$
  from which we see immediately that
  $\|S_\lambda\|=\|e_\lambda\|_p\|e_{-\lambda}\|_q$ and an easy
  calculation leads to the asymptotic formula given above.
  
  If $p>1$ then $f_\lambda$ is taken directly from the extremal case of 
  H\"older's inequality. If $p=1$ then there is no exact extremal function 
  for $S_\lambda$, but it is a simple calculation to check that
  $f_1$  is asymptotically extremal for any $g$ tending to $\infty$ at
  $\infty$.
\end{proof}

\begin{lem}\label{lem:ST}
  For any $p\in[1,\infty]$, the sequences of operators $(S_\lambda)$ and 
  $(T_\lambda)$ defined above are asymptotically equal as $\lambda\to\infty$ 
  through $\R^+$. In particular, $\|T_\lambda\|_p\sim C_p\e^{\lambda}/\lambda$.
\end{lem}
\begin{proof}
  Intuitively, $S_\lambda$ and $T_\lambda$ are close to each other for large
  $\lambda$ because their kernels differ only in the region $s>t$, where
  $\e^{\lambda(t-s)}$ is small when $\lambda$ is large. We can estimate the 
  rate of decay of $\|S_\lambda-T_\lambda\|$ as follows:
  \begin{align*}
    ((S_\lambda-T_\lambda)f)(t) &= \int_t^1 \e^{\lambda(t-s)}f(s)\ds \\
                                &= \int_0^{1-t}\e^{\lambda(t-1+u)}f(1-u)\du \\
                                &= \int_0^{1-t}\e^{-\lambda(1-t-u)}f(1-u)\du \\
                                &= (e_{-\lambda}*Rf)(1-t)
  \end{align*}
  where $R$ is the operator on $L^p(0,1)$ defined by $(Rf)(t)=f(1-t)$.
  We thus have that
  $$S_\lambda-T_\lambda=RV_{e_{-\lambda}}R.$$
  Now, $R$ is an isometric bijection on $L^p(0,1)$, so 
  $\|S_n-T_n\|_p=\|V_{e_{-\lambda}}\|_p$. We can now use the standard estimate
  to see that
  $$\|V_{e_{-\lambda}}\|_p\leq\int_0^1\e^{-\lambda t}\dt
    =\frac{1-\e^{-\lambda}}{\lambda}\sim\frac{1}{\lambda}$$
  so by Lemma~\ref{lem:rank-1}, 
  $\|S_\lambda-T_\lambda\|/\|S_\lambda\|\lesssim C_p^{-1}\e^{-\lambda}$ as 
  $\lambda\to\infty$. This shows that $S_\lambda\sim T_\lambda$ as
  $\lambda\to\infty$, so 
  $\|T_\lambda\|_p\sim\|S_\lambda\|_p\sim C_p\e^{\lambda}/\lambda$ by
  Lemma~\ref{lem:rank-1}.
\end{proof}

\begin{lem}\label{lem:exponential}
  For some fixed $\mu\in\R$ and $p\in[1,\infty]$, let $k(t)=\e^{\mu t}$ and 
  consider the Volterra operator $V_k$ acting on $L^p(0,1)$. Then 
  $$V_k^n\sim\frac{\e^{-(n-1)}}{\Gamma(n)}S_{n-1+\mu}$$
  and in particular
  $$\|V_k^n\|_p\sim\frac{C_p\e^\mu}{\Gamma(n+1)}$$
  as $n\to\infty$ through $\R^+$.
\end{lem}
\begin{proof}
  We have
  \begin{align*}
    \left\|\Gamma(n)V_k^n-\e^{-(n-1)}T_{n-1+\mu}\right\|_p 
      &\leq \int_0^1|t^{n-1}\e^{\mu t}-\e^{-(n-1)}\e^{(n-1+\mu)t}|\dt \\
      &=    \int_0^1 \e^{-(n-1)}\e^{(n-1+\mu)t}-t^{n-1}\e^{\mu t}\dt \\
      \intertext{(since $t^{n-1}\e^{\mu t}\leq\e^{-(n-1)}\e^{(n-1+\mu)t}$ for $t\in[0,1]$)}
      &\leq \e^\mu\int_0^1\e^{(n-1)(t-1)}-t^{n-1}\dt \\
      &<    \e^\mu\left(\frac{1}{n-1}-\frac{1}{n}\right) \\
      &=    \frac{\e^\mu}{n(n-1)}
  \end{align*}
  $$\frac{\|\Gamma(n)V_k^n-\e^{-(n-1)}T_{n-1+\mu}\|_p}
         {\|\e^{-(n-1)}T_{n-1+\mu}\|_p}\lesssim
    \frac{\e^\mu/n(n-1)}{\e^{-(n-1)}C_p\e^{n-1+\mu}/n}=
    \frac{1}{C_p(n-1)}\to 0.$$
  using Lemma~\ref{lem:ST}. This shows that
  $\Gamma(n)V_k^n\sim\e^{-(n-1)}T_{n-1+\mu}$ as $n\to\infty$. But
  $S_\lambda\sim T_\lambda$ as $\lambda\to\infty$ by Lemma~\ref{lem:ST},
  so 
  $$V_k^n\sim\frac{\e^{-(n-1)}}{\Gamma(n)}S_{n-1+\mu}$$
  as claimed. The asymptotic formula for $\|V_k^n\|$ now follows
  from Lemma~\ref{lem:parallel} and Lemma~\ref{lem:rank-1}.
\end{proof}

\begin{cor}\label{cor:power-exponential}
  Fix $\mu,r\in\R$ with $r>-1$, let $k(t)=t^r\e^{\mu t}$ and consider the 
  Volterra operator $V_k$ acting on $L^p(0,1)$ where $1\leq p\leq\infty$. Then 
  $$V_k^n\sim
    \frac{\Gamma(r+1)^n\e^{-((r+1)n-1)}}{\Gamma((r+1)n)}S_{(r+1)n-1+\mu}$$
  and in particular
  $$\|V_k^n\|_p\sim\frac{C_p\e^{\mu}(\Gamma(r+1))^n}{\Gamma((r+1)n+1)}$$
  as $n\to\infty$ through $\R^+$.
\end{cor}
\begin{proof}
  As remarked at the beginning of the section, if we define $k_0(t)=\e^{\mu t}$
  then we have $k^{*n}=(\Gamma(r+1))^nk_0^{*(r+1)n}$. The result is now
  immediate from Lemma~\ref{lem:exponential}.
\end{proof}

\section{More general kernels}

It is easy to see that if $h,k\in L^1(0,1)$ with $0\leq h\leq k$ then 
$\|V_h\|_p\leq\|V_k\|_p$ for any $p\in[1,\infty]$. This simple fact, in
combination with the results from  the previous section, allows us to
deduce asymptotic results for a large class of kernels.

\begin{lem}\label{lem:decay}
  Suppose $k$ is a measurable function on $[0,1]$ and there exist
  real constants $c,\mu,\nu,r$ with $c>0$ and $r>1$ such that
  $$ct^r\e^{\mu t}\leq k(t)\leq ct^r\e^{\nu t}$$
  for $t\in [0,1]$. Then for any $\delta\in(0,1)$, any $j\in\N$ and any 
  polynomial $P$,
  $$\frac{P(n)\int_0^{1-\delta}k^{*(n-j)}}{\|V_k^n\|_p}\to 0$$
  as $n\to\infty$.
\end{lem}
\begin{proof}
  Taking the $(n-j)$-fold convolution power of the right-hand inequality gives
  $$k^{*{n-j}}(t)\leq\frac{(c\Gamma(r+1))^n}{\Gamma((r+1)(n-j))}t^{(r+1)(n-j)-1}\e^{\nu t}$$
  Using the estimate $\e^{\nu t}\leq\max(1,\e^\nu)$, we have if $(r+1)(n-j)>1$,
  \begin{align*}
    \int_0^{1-\delta}k^{*(n-j)}
      &=    \frac{(c\Gamma(r+1))^{n-j}}{\Gamma((r+1)(n-j))}
            \int_0^{1-\delta}t^{(r+1)(n-j)-1}\e^{\nu t}\dt \\
      &\leq \frac{(c\Gamma(r+1))^{n-j}}{\Gamma((r+1)(n-j))}
            \frac{(1-\delta)^{(r+1)(n-j)}\max(1,\e^\nu)}{(r+1)(n-j)}.
  \end{align*}
  We can also see from the $n$-fold convolution power of the left-hand 
  inequality and Corollary~\ref{cor:power-exponential} that
  $$\|V_k^n\|_p\gtrsim\frac{C_p\e^{\mu}(c\Gamma(r+1))^n}{\Gamma((r+1)n+1)}.$$
  Combining these gives
  $$\frac{P(n)\int_0^{1-\delta}k^{*(n-j)}}{\|V_k^n\|_p}
    \lesssim\frac{P(n)\max(1,\e^\nu)\Gamma((r+1)n+1)(1-\delta)^{(r+1)(n-j)}}
                 {C_p\e^\mu(c\Gamma(r+1))^j\Gamma((r+1)(n-j))(r+1)(n-j)}.$$
  It is an immediate consequence of Stirling's formula that
  $\Gamma(n+s)/\Gamma(n)\sim n^s$ as $n\to\infty$ for any $s$, from
  which it follows that the right-hand side tends to zero as $n\to\infty$.
\end{proof}

We can now establish a localisation result: for a wide range of kernels, 
the asymptotic behaviour of $V_k^n$ is determined by the values of
$k$ in any neighbourhood of $0$.

\begin{lem}\label{lem:localisation}
  Suppose $h,k\in L^1(0,1)$, that $h$ and $k$ are equal on the interval
  $[0,\delta]$ for some $\delta\in(0,1)$ and that there exist real constants 
  $c,\mu,\nu,r$ with $c>0$ and $r>1$ such that
  $$ct^r\e^{\mu t}\leq h(t)\leq ct^r\e^{\nu t}$$
  for $t\in [0,1]$. Then for any $p\in[1,\infty]$, $V_k^n\sim V_h^n$ on
  $L^p(0,1)$ as $n\to\infty$.
\end{lem}
\begin{proof}
  Let $g=k-h$, so $k=h+g$ and $g$ is zero on $[0,\delta]$.
  We can use the binomial theorem in the convolution algebra $L^1(0,1)$ 
  to write
  $$k^{*n}=(h+g)^{*n}=
    h^{*n}+g^{*n}+\sum_{j=1}^{n-1}\binom{n}{j}g^{*j}*h^{*(n-j)}.$$
  Now, if we were working on the whole of $\R$, then $g^{*n}$ 
  would be supported on $[n\delta,n]$ and $g^{*j}*h^{*(n-j)}$ would be 
  supported on $[j\delta,n]$. But we are working in $L^1(0,1)$, so if we 
  choose $N>1/\delta$ then for $n\geq N$ we have
  $$k^{*n}=h^{*n}+
    \sum_{j=1}^{N-1}\binom{n}{j}g^{*j}*h^{*(n-j)}.$$
  Moreover, since $g^{*j}$ is supported to the right of $j\delta$, we have
  $$g^{*j}*h^{*(n-j)}=
    g^{*j}*(h^{*(n-j)}\chi_{[0,1-j\delta]})$$
  and hence
  $$k^{*n}-h^{*n}=\sum_{j=1}^{N-1}\binom{n}{j}
    g^{*j}*(h^{*(n-j)}\chi_{[0,1-j\delta]}).$$
  We can therefore estimate
  $$\frac{\|V_k^n-V_h^n\|_p}{\|V_h^n\|_p}\leq\sum_{j=1}^N\binom{n}{j}
    \frac{\left(\int_0^1|g^{*j}|\right)\left(\int_0^{1-j\delta}h^{*(n-j)}\right)}
         {\|V_h^n\|_p}.$$
  This is a finite sum of terms, all of which tend to zero by
  Lemma~\ref{lem:decay}, so we can conclude that 
  $V_h^n\sim V_k^n$ as $n\to\infty$.
\end{proof}

We are now in a position to prove the main result.

\begin{thm}
\label{thm:main}
  Suppose $k\in L^1(0,1)$ is such that $k(t)=t^rf(t)$ where $r>-1$,
  $f(0)\neq 0$ and $f'(0)$ exists, and let
  $$h(t)=f(0)t^r\e^{(f'(0)/f(0))t}.$$
  Then for any $p\in[1,\infty]$, $V_k^n\sim V_h^n$ on $L^p(0,1)$. It
  follows that $V_k^n$ is also asymptotically equivalent to the
  sequence of rank $1$ operators described in 
  Corollary~\ref{cor:power-exponential}; in particular:
  $$\|V_k^n\|_p\sim
    \frac{C_p(|f(0)|\Gamma(r+1))^n\e^{f'(0)/f(0)}}{\Gamma((r+1)n+1)}$$
  and if 
  $$f_n(t)=\begin{cases}
             \e^{-g(n)nt}                       & \text{ if } p=1 \\
             \e^{-((r+1)n-1+k'(0)/k(0))t/(p-1)} & \text{ if } 1<p<\infty \\
             1                                  & \text{ if } p=\infty
           \end{cases}$$
  where $g$ is any function such that $g(n)\to\infty$ as 
  $n\to\infty$, then $(f_n)$ is asymptotically extremal for $(V_k^n)$.
\end{thm}
\begin{proof}
  For $\eta\in\R$, let
  $$h_\eta(t)=f(0)t^r\e^{(f'(0)/f(0)+\eta)t}.$$
  We can assume without loss of generality that $f(0)>0$, so $\log f$
  is differentiable at $0$ and hence for any $\eta>0$ there exists
  $\delta_\eta\in(0,1)$ such that if $0<t\leq\delta_\eta$ then
  $$(\log f)'(0)-\eta\leq\frac{\log f(t)-\log f(0)}{t}\leq (\log f)'(0)+\eta$$
  or equivalently
  $$f(0)t^r\e^{(f'(0)/f(0)-\eta)t}\leq k(t)\leq f(0)t^r\e^{(f'(0)/f(0)+\eta)t}.$$
  Now let
  $$k_\eta(t)=\begin{cases}
                k(t) & \text {if } 0\leq t\leq \delta_\eta \\
                h(t) & \text {if } \delta_\eta<t\leq 1
              \end{cases}$$
  so $h_{-\eta}\leq k_\eta\leq h_\eta$ and $h_{-\eta}\leq h\leq h_\eta$.
  Because all the functions involved are non-negative, we can take the
  $n$-fold convolution power of these inequalities to give 
  $h_{-\eta}^{*n}\leq k_\eta^{*n}\leq h_\eta^{*n}$ and 
  $h_{-\eta}^{*n}\leq h^{*n}\leq h_\eta^{*n}$. If follows that
  $|k_\eta^{*n}-h^{*n}|\leq h_{\eta}^{*n}-h_{-\eta}^{*n}$ and we can
  integrate to give, abbreviating $f'(0)/f(0)$ to $\mu$,
  \begin{align*}
    \|V_{k_\eta}^n-V_h^n\|_p 
      &\leq \frac{(f(0)\Gamma(r+1))^n}{\Gamma((r+1)n)}
            \int_0^1 t^{(r+1)n-1}\left(\e^{(\mu+\eta)t}-\e^{(\mu-\eta)t}\right)\dt \\
      &\leq \frac{(f(0)\Gamma(r+1))^n}{\Gamma((r+1)n)}
            \left(\e^{\mu+\eta}-\e^{\mu-\eta}\right)\int_0^1 t^{(r+1)n-1}\dt \\
      &=    \frac{(f(0)\Gamma(r+1))^n}{\Gamma((r+1)n)}
            \frac{1}{(r+1)n}\left(\e^{\mu+\eta}-\e^{\mu-\eta}\right) \\
      &\leq \frac{K_1(f(0)\Gamma(r+1))^n\eta}{\Gamma((r+1)n+1}
  \end{align*}
  for all $n\in\N$ and all $\eta\in[0,1]$, say, where $K_1$ is a constant
  independent of $n$ and $\eta$. We also have
  $$k_\eta^{*n}\geq h_{-\eta}^{*n}\geq\e^{-\eta}h^{*n}$$
  so $\|V_{k_\eta}^n\|_p\geq\e^{-\eta}\|V_h^n\|_p$. But 
  $$\|V_h^n\|_p\sim\frac{C_p(f(0)\Gamma(r+1))^n\e^\mu}{\Gamma((r+1)n+1)}$$
  so in particular
  $$\|V_h^n\|_p\geq\frac{K_2(f(0)\Gamma(r+1))^n\e^\mu}{\Gamma((r+1)n+1)}$$
  for all $n\in\N$, where $K_2$ is independent of $n$. Combining all these,
  we see that
  $$\frac{\|V_{k_\eta}^n-V_h^n\|_p}{\|V_{k_\eta}^n\|_p}\leq K_3\eta\e^{\eta}$$
  for all $n\in\N$ and all $\eta\in[0,1]$, where $K_3$ is independent of
  $n$ and $\eta$.
  
  Now, for any $\eps>0$ we can find $\eta\in(0,1)$ such that
  $$\frac{\|V_{k_\eta}^n-V_h^n\|}{\|V_{k_\eta}^n\|}<\frac{\eps}{2\e}$$
  for all $n\in\N$. We can also use Lemma~\ref{lem:localisation} to find 
  $N\in\N$ such that if $n>N$ then
  $$\frac{\|V_k^n-V_{k_\eta}^n\|_p}{\|V_{k_\eta}^n\|_p}<\frac{\eps}{2\e}$$
  and hence
  $$\frac{\|V_k^n-V_h^n\|_p}{\|V_{k_\eta}^n\|_p}<\frac{\eps}{\e}.$$
  But $k_\eta^{*n}\leq h_{\eta}^{*n}\leq\e^\eta h^{*n}\leq \e h^{*n}$
  since $\eta\in(0,1)$. We therefore have $\|V_{k_\eta}^n\|\leq
  \e\|V_h^n\|$, so for $n>N$ we have
  $$\frac{\|V_k^n-V_h^n\|_p}{\|V_h^n\|_p}<\eps$$
  showing that $(V_k^n)$ and $(V_h^n)$ are asymptotically equal. Their
  norms are thus also asymptotically equal so we have
  $$\|V_k^n\|_p\sim\frac{C_p(f(0)\Gamma(r+1))^n\e^\mu}{\Gamma((r+1)n+1)}$$
  by Corollary~\ref{cor:power-exponential}, as claimed. We also know from
  Corollary~\ref{cor:power-exponential} that
  $$V_h^n\sim
    \frac{\Gamma(r+1)^n\e^{-((r+1)n-1)}}{\Gamma((r+1)n)}S_{(r+1)n-1+\mu}$$
  as $n\to\infty$, where $S_\lambda$ and $T_\lambda$ is as defined
  in Section~\ref{sec:power-exponential}, so we have
  $$V_k^n\sim 
    \frac{\Gamma(r+1)^n\e^{-((r+1)n-1)}}{\Gamma((r+1)n)}S_{(r+1)n-1+\mu}.$$
  By Lemma~\ref{lem:parallel}, these two sequences of operators have the same
  asymptotically extremal sequences of vectors. An appropriate sequence
  for $(S_\lambda)$ was identified in Lemma~\ref{lem:rank-1};
  substituting $\lambda=(r+1)n-1+k'(0)/k(0)$ gives the sequence in the
  statement of the theorem.
\end{proof}

\section{Further remarks on the case $p=1$: a probabilistic interpretation}
\label{sec:large}

In the case $p=1$, the estimate used throughout is in fact exact:
$\|V_k\|_1=\|k\|_1$ (consider the action of $V_k$ on an approximate identity).
Theorem~\ref{thm:main} thus gives the following result about powers of
elements of the Volterra algebra $L^1(0,1)$:

\begin{cor}
  Suppose $k\in L^1(0,1)$ is such that $k(t)=t^rf(t)$ where $r>-1$,
  $f(0)\neq 0$ and $f'(0)$ exists. Then
  $$\|k^{*n}\|_1\sim
    \frac{(|f(0)|\Gamma(r+1))^n\e^{f'(0)/f(0)}}{\Gamma((r+1)n+1)}$$
  as $n\to\infty$.
\end{cor}

If $k\in L^1(0,\infty)$, $k\geq 0$ a.e.\ and $\int_0^\infty k=1$ then we
can interpret $k$ as the probability density of a random variable and
$k^{*n}$ as the density of the sum of $n$ independent random variables with 
density $k$. The $L^1$ norm of the restriction to $(0,1)$ of $k^{*n}$ is
then the probability that this sum is no larger than $1$.

\begin{cor}
  Suppose $k\in L^1(0,\infty)$ is a probability density function and that 
  $k(t)=t^rf(t)$ where $r>-1$, $f(0)\neq 0$ and $f'(0)$ exists. Let
  $(X_n)$ be a sequence of independent random variables with this
  density, and let $S_n=X_1+X_2+\dots+X_n$. Then
  $$\mathbf{P}(S_n\leq 1)\sim
    \frac{(f(0)\Gamma(r+1))^n\e^{f'(0)/f(0)}}{\Gamma((r+1)n+1)}$$
  as $n\to\infty$.
\end{cor}

This limit theorem seems to go beyond the scope of known results on
such sums, such as those in Petrov~\cite[Section~5.8]{Petrov95}. In the
notation of that section, we have $x=O(n^{1/2})$ but not $x=o(n^{1/2})$ 
which, as explicitly noted, is not sufficient for the results there to apply.

\end{document}